\documentclass[11pt]{article}
\usepackage{times}
\usepackage{fullpage}
\usepackage{amsmath}
\usepackage{amssymb}
\usepackage{amsfonts}
\usepackage{epsfig}
\usepackage[dvips]{color}
\usepackage{colortbl}

\newcounter{algoctr}

\newif\ifnotesw\noteswtrue
   {\ifnotesw\marginpar[\hfill\(\top\)]{\(\top\)}\fi}%
   {\ifnotesw\marginpar[\hfill\(\bot\)]{\(\bot\)}\fi}

\newcommand{\mnote}[1]%
    {\ifnotesw\marginpar%
        [{\scriptsize\begin{minipage}[t]{\marginparwidth}
        \raggedleft#1%
                        \end{minipage}}]%
        {\scriptsize\begin{minipage}[t]{\marginparwidth}
        \raggedright#1%
                        \end{minipage}}%
    \fi}

\newcommand{\ignore}[1]{}

\newcommand{\etal}{{\it et al.~}}

\newsavebox{\given}
\savebox{\given}[1em]{\rule[-1.5ex]{.2mm}{4ex}}

\newtheorem{theorem}{Theorem}

\newtheorem{lemma}[theorem]{Lemma}
\newtheorem{fact}[theorem]{Fact}

\newcommand{\blackslug}{\rule{7pt}{7pt}}

\newcommand{\qed}{\hfill{\setlength{\fboxsep}{0pt}
\framebox[7pt]{\rule{0pt}{7pt}}}}

\renewcommand{\notin}{\ifmmode \not\in \else $\not\in$ \fi}
\newlength{\thislabel}
\newcommand{\labsize}[1]{\settowidth{\thislabel}{#1}}

\newcommand{\prf}{\par\noindent{\sl Proof } \hspace{.01 in}}

\newcommand{\bra}[1]{\langle #1 |}
\newcommand{\ket}[1]{| #1 \rangle}
\newcommand{\braket}[2]{\langle #1 | #2 \rangle}

\newcommand{\Real}{\mathbb{R}}

\newcommand{\va}{\mathbf{a}}
\newcommand{\vb}{\mathbf{b}}

\newcommand{\vu}{\mathbf{u}}
\newcommand{\vv}{\mathbf{v}}
\newcommand{\ta}{\mathcal{A}}
\newcommand{\tb}{\mathcal{B}}
\newcommand{\tc}{\mathcal{C}}
\newcommand{\td}{\mathcal{D}}

\newcommand{\ttt}{\mathcal{T}}

\newcommand{\tx}{\mathcal{X}}
\newcommand{\tl}{\mathcal{L}}

\newcommand{\GG}{\mathcal{G}}

\newcommand{\Ring}{\mathsf{R}}

\newcommand{\Vector}{\mathbb{V}}
\newcommand{\Matrix}{\mathbb{M}}
\newcommand{\Module}{\mathfrak{M}}
\newcommand{\Tensor}{\mathbb{T}}
\newcommand{\Zn}{\mathbb{Z}_{n}}

\newcommand{\scalarprod}{\circ}

\newcommand{\basis}{\mathfrak{B}}

\title{
Tensors as module homomorphisms over group rings
} 
\author{
Carmeliza Navasca\footnote{Department of Mathematics, Clarkson University. email: cnavasca@clarkson.edu}
\and
Michael Opperman\footnote{Department of Mathematics, Clarkson University. email: oppermmc@clarkson.edu}
\and
Timothy Penderghest\footnote{Department of Computer Science, Clarkson University. email: pendertj@clarkson.edu}
\and
Christino Tamon\footnote{Department of Computer Science, Clarkson University. email: tino@clarkson.edu. Contact author.}
}
\date{\today}

\begin{document}
\bibliographystyle{plain}
\maketitle

\begin{abstract}
Braman \cite{b08} described a construction where third-order tensors are exactly the set of 
linear transformations acting on the set of matrices with vectors as scalars. This extends the familiar 
notion that matrices form the set of all linear transformations over vectors with real-valued scalars. 
This result is based upon a circulant-based tensor multiplication due to Kilmer \etal \cite{kmp08}.
In this work, we generalize these observations further by viewing this construction in its natural 
framework of group rings. 
The circulant-based products arise as convolutions in these algebraic structures. 
Our generalization allows for any abelian group to replace the cyclic group, 
any commutative ring with identity to replace the field of real numbers,
and an arbitrary order tensor to replace third-order tensors, provided the underlying ring is commutative.

\vspace{0.05in}
\end{abstract}

\section{Introduction}

Matrix multiplication is an example of a {\em contracted} product between two second-order tensors.
Namely, for two matrices $A = (a_{ij})$ and $B = (b_{jk})$, where $A,B \in \Real^{n \times n}$, 
their product $C = AB$ is defined as 
\begin{equation}
c_{ik} = \sum_{j} a_{ij}b_{jk},
\ \ 
\mbox{ or simply }
\ \ 
c_{ik} = a_{ij}b_{jk}.
\end{equation}
The latter notation was attributed to Einstein, where a summation is implied whenever two indices appear exactly twice.
With this notation, the inner product of two vectors $\vu = (u_{j})$ and $\vv = (v_{j})$ is denoted by $u_{j}v_{j}$
and the matrix-vector product of $A$ and $\vu$ is denoted by $a_{ij}u_{j}$.
It is possible to extend this contracted product to higher-order tensors, but it lacks the closure property 
on odd-order tensors (as the example on vectors showed). 
So, it is not possible to define a tensor multiplication for odd-order tensors based on contracted product 
if we require the closure property.

Recently, Kilmer \etal \cite{kmp08} proposed an interesting tensor multiplication for third-order tensors
based on circulant matrices. A standard $\Zn$-circulant $n \times n$ matrix $A$ (see \cite{d94})
is defined completely by its first row as $a_{i,j} = a_{0,j-i} = a_{i-j}$, where subtraction is done modulo $n$: 
\begin{equation}
A = circ(\langle a_{0},a_{1},\ldots,a_{n-1}\rangle) = 
\begin{bmatrix}
a_{0}   & a_{n-1}  & a_{n-2}  & \ldots & a_{1} \\
a_{1} & a_{0}  & a_{n-1}  & \ldots & a_{2} \\
\vdots  & \vdots & \vdots & \vdots & \vdots  \\
a_{n-1}   & a_{n-2}  & a_{n-3}  & \ldots & a_{0}
\end{bmatrix}.
\end{equation}
We view third-order tensors $\ta,\tb \in \Real^{n \times n \times n}$ as a sequence of matrices 
$\ta = (A_{k})$ and $\tb = (B_{k})$, where $A_{k},B_{k} \in \Real^{n \times n}$. 
Then, their {\em circulant} product $\tc = \ta \star \tb$ is given by
\begin{equation}
\begin{bmatrix}
A_{0}   & A_{n-1}  & A_{n-2}  & \ldots & A_{1} \\
A_{1} & A_{0}  & A_{n-1}  & \ldots & A_{2} \\
\vdots  & \vdots & \vdots & \vdots & \vdots  \\
A_{n-1}   & A_{n-2}  & A_{n-3}  & \ldots & A_{0}
\end{bmatrix}
\begin{bmatrix}
B_{0}  \\
B_{1}  \\
\vdots \\
B_{n-1}
\end{bmatrix}
=
\begin{bmatrix}
C_{0}  \\
C_{1}  \\
\vdots \\
C_{n-1}
\end{bmatrix}
\end{equation}
where $\tc = (C_{k})$ (viewed as a sequence of matrices) with $C_{k} \in \Real^{n \times n}$. 
This tensor multiplication has the important closure property lacking in a contracted product 
and allows a singular value decomposition to be defined easily for third-order tensors (see \cite{kmp08}).
We note that $\tc = (C_{k})$ is simply the {\em convolution} between the two matrix sequences
$\ta = (A_{k})$ and $\tb = (B_{k})$, since $C_{k} = \sum_{r+s=k} A_{r}B_{s}$. 
By viewing a matrix as a sequence of vectors, we can use the circulant-based tensor multiplication 
to define a tensor-matrix product similarly.

A fundamental fact in linear algebra states that the set of $n \times n$ matrices is isomorphic to
the set of linear transformations over the vector space $\Real^{n}$. A matrix $A$ linearly maps a vector 
$\vu$ to another vector $\vv$ through a matrix-vector product $\vv = A\vu$ ~or~ $v_{i} = a_{ij}u_{j}$, 
defined using a contracted product. This observation is difficult to generalize for odd-order tensors 
since the matrix-vector product does not map a vector space to itself.
Building on the work in \cite{kmp08}, Braman \cite{b08} proved the following surprising property
of third-order tensors: the set of all third-order tensors form precisely the set of linear transformations on 
a module (vector space over ring) defined by second-order tensors (matrices) with scalars from first-order tensors
(vectors). In this sense, Braman's result recovered the fundamental isomorphism between matrices and linear 
transformations on a vector space (over a field). This answered one of the main questions posed in \cite{martin}:
{\em to characterize exactly the objects operated on by (third-order) tensors}.

Braman's construction relied on a set of carefully defined products between vectors (namely, $\va \odot \vb$),
between vectors and matrices (namely, $\va \scalarprod X$), and also the circulant product between third-order
tensors (namely, $\ta \star \tb$) which also includes a product between third-order tensors and matrices (namely, $\ta \star X$).
Figure \ref{figure:products} contains a description of these products and their definitions in \cite{b08}.
Using these products, Braman proved that the set of matrices forms a unitary free module with scalars that are vectors. 
Moreover, any third-order tensor is a linear transformation over this module and any linear transformation over the
module can be represented as a third-order tensor. An interesting diagonalization theorem for third-order tensors 
was also proved which recovers the eigenvalue-eigenvector equations for diagonalizable matrices.
This answered another main question posed in \cite{martin}: {\em to determine the singular or eigenvalues of tensors 
and to see if these are possibly vectors rather than scalars}.

\begin{figure}[t]
\begin{center}
\begin{tabular}{|c||c||c|c|} \hline
product		&	$\Zn$-circulant matrix	&	abstract convolution 	&	relevant group ring \\	\hline
$\va \odot \vb$ 	& 	$circ(\va) \cdot \vb$ 	&	$\va \star \vb$			&	$\Vector = \Ring\GG$	\\
\rowcolor[gray]{0.85}
$\va \scalarprod X$	& 	$X \cdot circ(\va)$ 	&	$\va \star X$			&	$\Ring\GG$ and $\Vector\GG$ \\
$\ta \star X$   	&	$circ(\ta) \cdot X$		&	$\ta \star X$			&	$\Tensor = \Matrix\GG$	\\	
$\ta \star \tb$   	&	$circ(\ta) \cdot \tb$	&	$\ta \star \tb$			&	$\Tensor = \Matrix\GG$	\\	\hline
\end{tabular}

\vspace{.1in}
{Notation}: $\va,\vb$ are ``vectors'', $X$ is a ``matrix'', $\ta,\tb$ are ``tensors''.
\caption{A comparison of tensor products defined by Braman \cite{b08}, which is based on
circulant matrices, and our interpretation, which is based on convolutions in group rings. 
The ingredients in our construction are:
an abelian group $\GG$,
a commutative ring $\Ring$ with identity,
a group ring $\Vector = \Ring\GG$,
a matrix ring $\Matrix = M_{n}(\Ring)$ (viewed also as a group ring $\Vector\GG$),
and a group ring $\Tensor = \Matrix\GG$.}
\label{figure:products}
\end{center}
\hrule
\end{figure}

In this work, we generalize Braman's result \cite{b08} in several ways.
First, we show that all of the different products defined in the construction are, in fact, convolutions 
in various group rings. This provides a unifying picture to Braman's construction on third-order tensors.
Then, we extend the construction to arbitrary order tensors (beyond third-order) by viewing the set of tensors 
equipped with addition and convolution product as a (commutative) ring with identity. This extends previous 
constructions which defined the tensors over the field of real or complex numbers.
Our generalization also allows the underlying circulant matrices to be generalized to any {\em abelian} 
group-theoretic circulant.
Within our algebraic setting, we also identified Braman's clever choice of a basis for the unitary free
module of matrices, which is not the natural basis for this group ring.
Figure \ref{figure:products} describes the new interpretation in our framework 
of the tensor products used in \cite{b08}.

This paper is organized as follows. In Section \ref{section:preliminaries}, we describe relevant notation and 
algebraic background which will be used throughout. 
In Section \ref{section:tower}, we provide a unifying view of Braman's construction by defining a tower of group rings
which forms the natural setting for her construction. Some additional products which arise naturally from the tower
of rings are described in Section \ref{subsection:induced-products}. 
In Section \ref{section:unitary-free-module}, we prove that the set 
of matrices is a unitary free module under a special scalar product between vectors and matrices.
In Section \ref{section:isomorphism}, we show that the set of tensors is isomorphic to the set of module homomorphisms 
over the set of matrices. Finally, in Section \ref{section:diagonalization}, we show that the module structure on
tensors admits a particular form of tensor diagonalization.

\section{Preliminaries} \label{section:preliminaries}

\paragraph{Notation}
The Kronecker delta $\delta_{j,k}$ is defined to be $1$ if $j=k$ and $0$ otherwise.
We will adopt the convenient Dirac notation for matrices. Here, $\ket{k}$ will be used to denote 
a $n \times 1$ matrix (column vector) which is one at position $k$ and zero elsewhere. 
The corresponding $1 \times n$ matrix (row vector) will be denoted as $\bra{k}$. 
Thus, $\bra{j}X\ket{k}$ denotes the $(j,k)$-entry of a $n \times n$ matrix $X$. 
As is standard, the Dirac bracket (inner product) is $\braket{j}{k} = \delta_{j,k}$
while the outer product $\ket{j}\bra{k}$ is a matrix whose $(j,k)$-entry is $1$ and the other entries are zero.\\

For a group $\GG$ and a ring $\Ring$ with identity $1_{\Ring}$, the {\em group ring} $\Ring\GG$ 
is the set of all formal sums $\sum_{g \in \GG} a_{g} \ket{g}$, where $a_{g} \in \Ring$.
Here, we use the notation $\ket{g}$ to denote a placeholder for the group element $g$, 
which may also be viewed as the unit vector which is one only at position indexed by $g$ and zero elsewhere.
If $1_{\GG}$ is the identity of $\GG$ and $a \in \Ring$, then $a\ket{1_{\GG}}$ is simply written as $a$.
Also, $1_{\Ring}\ket{g}$, for $g \in \GG$, is simply denoted by $\ket{g}$. 
The group ring $\Ring\GG$ is a ring under component-wise addition and convolution-like multiplication.
More specifically, addition is defined as:
\begin{equation}
(\sum_{g \in \GG} a_{g} \ket{g}) \ + \ (\sum_{g \in \GG} b_{g} \ket{g})
    \ = \ \sum_{g \in \GG} (a_{g}+b_{g}) \ket{g},
\end{equation}
and multiplication is defined as:
\begin{equation}
(\sum_{g \in \GG} a_{g} \ket{g}) \ \star \ (\sum_{h \in \GG} b_{h} \ket{h}) 
    \ = \ \sum_{k \in \GG} c_{k} \ket{k}, \ \ \
	        \mbox{ where } c_{k} = \sum_{\substack{g,h \in \GG:\\ gh = k}} a_{g}b_{h},
\end{equation} 
where we let $(a\ket{g}) \star (b\ket{h}) = (ab)\ket{gh}$, for $a,b \in \Ring$ and $g,h \in \GG$. 
This multiplication is simply a {\em convolution} (with respect to the group $\GG$) between 
the ``sequences'' $a$ and $b$.
It is known that $\Ring\GG$ is a commutative ring with identity 
if and only if $\Ring$ is commutative with identity and $\GG$ is abelian (see \cite{h74}, page 117). 

Let $\Ring$ be a ring with identity. A {\em (left) $\Ring$-module} $\Module$ is an abelian group together 
with a scalar multiplication $\scalarprod : \Ring \times \Module \rightarrow \Module$ which satisfies
the following axioms. For any $\va,\vb \in \Ring$ and $X,Y \in \Module$, we have: 
\begin{eqnarray}
\va \scalarprod (X + Y) & = & \va \scalarprod X + \va \scalarprod Y \\
(\va + \vb) \scalarprod X & = & \va \scalarprod X + \vb \scalarprod X \\
(\va \star \vb) \scalarprod X & = & \va \scalarprod (\vb \scalarprod X).
\end{eqnarray}
The module $\Module$ is called {\em unitary} if there is an element $1_{\Ring} \in \Ring$ so that
$1_{\Ring} \scalarprod X = X$, for all $X \in \Module$. 
$\Module$ is called a {\em free} module if there is a linearly independent set 
$\mathfrak{B} = \{B_{1},\ldots,B_{r}\} \subseteq \Module$, so that any element $X \in \Module$ can be written as 
$\sum_{k=1}^{r} \va_{k} \scalarprod B_{k}$, for some $\va_{k} \in \Ring$. Such a set $\mathfrak{B}$ is called a
{\em basis} for $\Module$.

If $\Module_{1}$ and $\Module_{2}$ are modules over a ring $\Ring$, a function 
$\psi: \Module_{1} \rightarrow \Module_{2}$ is an
{\em $\Ring$-module homomorphism} if for all $X,Y \in \Module_{1}$ and $\va \in \Ring$:
\begin{eqnarray}
\psi(X + Y) & = & \psi(X) + \psi(Y) \\
\psi(\va \scalarprod X) & = & \va \circ \psi(X).
\end{eqnarray}
If $\Ring$ is a division ring, then $\psi$ is called a {\em linear transformation} (see \cite{h74}).
For a module $\Module$, we denote the set of all $\Ring$-module homomorphisms from $\Module$ 
to itself as $\hom_{\Ring}(\Module)$ (also called {\em endomorphisms}). 
Note $\hom_{\Ring}(\Module)$ is a ring under component-wise addition and composition of maps.

\begin{figure}[t]
\begin{center}
\epsfig{file=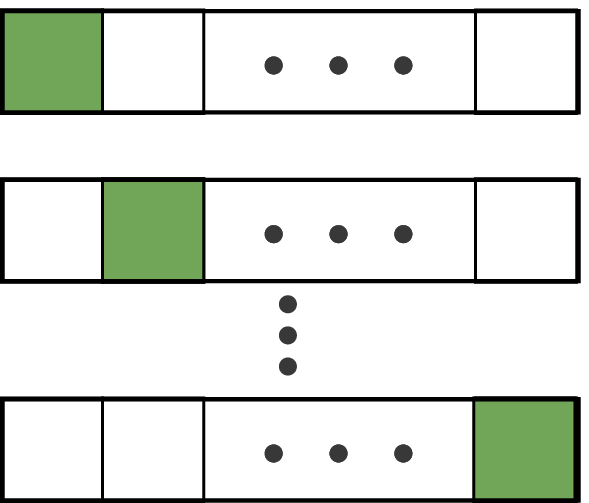, scale=0.5}
\caption{The group ring basis $\{\ket{g} : g \in \Zn\}$ of $\Vector = \Ring\Zn$.}
\label{figure:vector_basis}
\end{center}
\end{figure}

\begin{fact} (Ribenboim \cite{r69})
For a commutative ring $\Ring$ with identity, a group ring $\Ring\GG$ is unitary free module.
A natural basis for $\Ring\GG$ is $\{\ket{g} : g \in \GG\}$, where $\ket{g}$ is shorthand for $1_{\Ring}\ket{g}$.
\end{fact}

\section{Tower of rings} \label{section:tower}

\begin{figure}[t]
\begin{center}
\epsfig{file=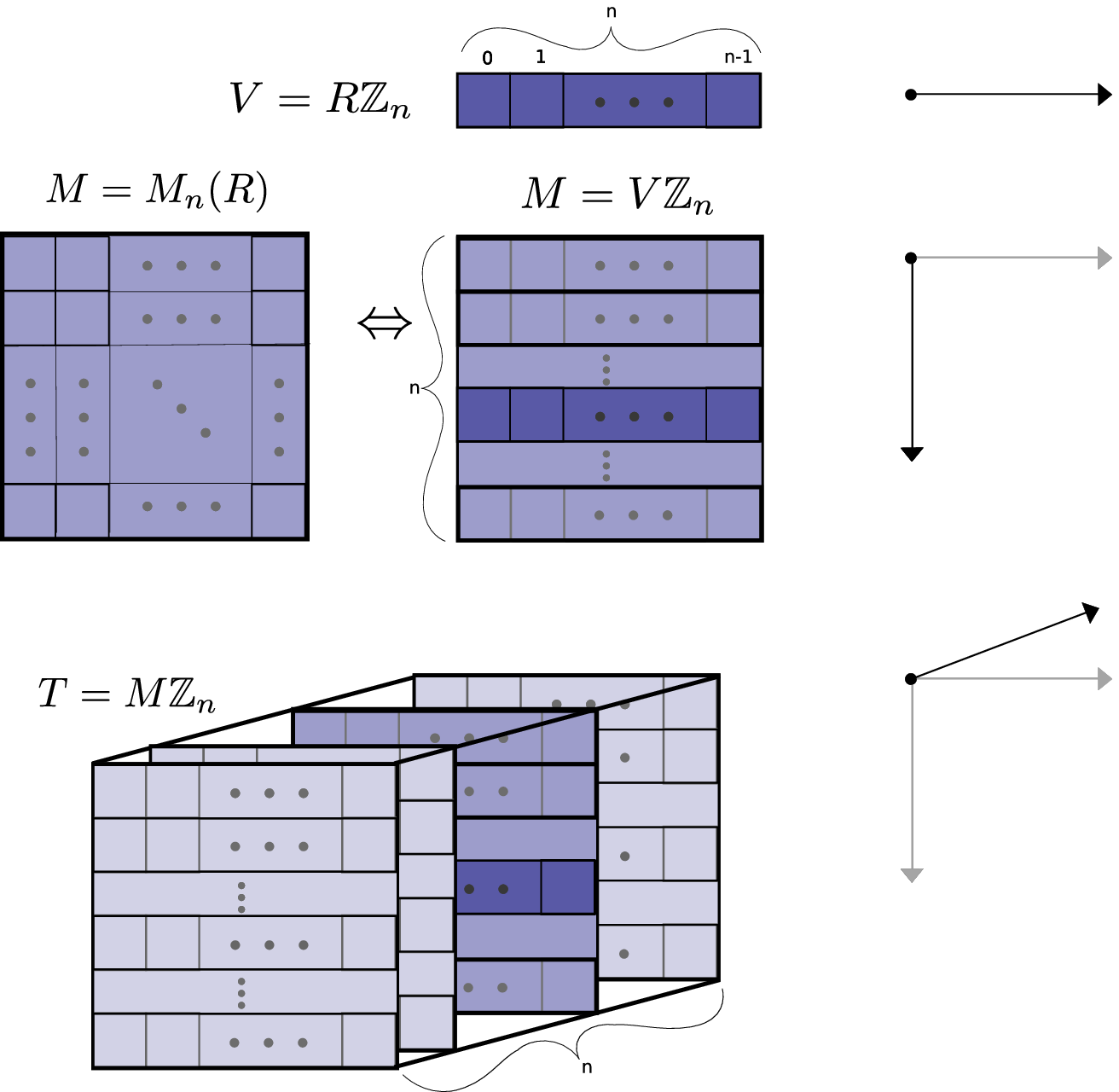, scale=0.75}
\caption{Tower of rings: the case when $\GG = \Zn$.}
\label{figure:matrices}
\end{center}
\end{figure}

We view Braman's construction \cite{b08} in an algebraic framework involving group rings.
\begin{enumerate}
\item 
Let $\GG$ be an abelian group $\GG$ of order $n$. \\
The choice of $\GG$ determines the type of convolution or circulant used.
A standard choice is the cyclic group $\Zn$ of order $n$ as used in \cite{kmp08,b08}.

\item
Let $\Ring$ be a commutative ring with identity. \\
This choice reflects the underlying algebraic ring structure and it generalizes the case where
$\Ring = \Real$ is the field of real numbers.

\item
Let $\Vector$ be the group ring $\Ring\GG$. \\
This is a commutative ring with identity whenever $\GG$ is an abelian group 
and $\Ring$ is a commutative ring with identity (see \cite{h74}, page 117).
A natural example here is $\Vector = \Real^{n}$ which is the $n$-dimensional vector space over the reals.

\item
Let $\Matrix$ be the matrix ring $M_{n}(\Ring)$. \\
This ring admits two alternate interpretations:
\begin{enumerate}
\item
$\Matrix \cong \hom_{\Ring}(\Vector)$: 
$\Matrix$ is the ring of all module-homomorphisms over $\Vector$ (see \cite{h74}, p330).

\item 
$\Matrix \cong \Vector\GG$:
$\Matrix$ is a group ring built from the ring $\Vector$ and the abelian group $\GG$.
\end{enumerate}
Note the matrix ring and the group ring interpretations are not algebraically equivalent since 
the former is not a commutative ring while the latter is. 
A standard example here is $\Matrix = M_{n}(\Real)$ or the set of $n \times n$ real matrices.

\item
Let $\Tensor$ be the group ring $\Matrix\GG \cong \hom_{\Ring}(\Vector)\GG$. \\
A main example is the set of $n \times n \times n$ tensors over the reals.

\end{enumerate}
We will denote a convolution in a group ring generically by $\star$.
Note this generic construction can be extended to arbitrary order tensors
by letting $\Ring$ be a {\em commutative} ring of tensors with identity.

\subsection{Induced products} \label{subsection:induced-products}

We describe some additional products which are induced from the rings involved.
\begin{itemize}
\item
The natural product among tensors in $\Tensor = \hom_{\Ring}(\Vector)\GG$ is simply their 
(convolution) product in the group ring itself. This product recursively uses the composition map 
between module-homomorphisms in $\hom_{\Ring}(\Vector)$.

\item
A natural tensor-matrix product between elements of $\Tensor$ and $\Matrix$ can be defined as follows.
Let $\star: \Tensor \times \Matrix \rightarrow \Matrix$ be a convolution between 
an element of $\Tensor = \hom_{\Ring}(\Vector)\GG$ and an element of $\Matrix = \Vector\GG$. 
This product recursively uses the module-homomorphism action on $\Vector$ in a natural way.
Alternatively, we may view this tensor-matrix product as the (convolution) product in $\Tensor$ 
since $\Vector$ is isomorphic to a collection of ``constant'' homomorphisms
$\{\phi_{\vv} : \vv \in \Vector\}$, where $\phi_{\vv}(x) = \vv$, for all $x \in \Vector$.

\item
A ``scalar'' product between $\va = \sum_{r} \va_{r}\ket{r} \in \Vector$ 
and $X = \sum_{s} X_{s}\ket{s} \in \Matrix$ is defined as a mixed convolution $\va \scalarprod X = \va \star X$, where
\begin{equation}
\va \star X = \sum_{g} (\sum_{rs=g} \va_{r} \star X_{s}) \ket{g},
\end{equation}
and the product $\va_{r} \star X_{s}$ is a convolution in $\Vector$ 
(viewing $\va_{r} \in \Ring$ as $\va_{r}\ket{1_{\GG}} \in \Vector$).

\end{itemize}
{\em Remark}: Braman \cite{b08} used a ring anti-isomorphism $\phi$ of $\Vector$ (see \cite{h74}, page 330)
to define the scalar product as $\va \scalarprod X = \phi(\va) \star X$, where 
$\phi(\sum_{g} \va_{g}\ket{g}) = \sum_{g} \va_{g}\ket{g^{-1}}$.
This definition is equivalent to~ $\va \scalarprod X = X \cdot circ(\va)$ and 
reflects the transposition present in the matrix multiplication.

\section{Unitary free module} \label{section:unitary-free-module}

We prove some properties showing the interplay between the distinct products involving 
vectors, matrices and tensors, and show that $\Matrix$ has a vector space structure over
the underlying ring.

\begin{lemma} \label{lemma:scalar-assoc}
For any $\va,\vb \in \Vector$ and any $X \in \Matrix$, we have
\begin{equation}
\va \scalarprod (\vb \scalarprod X) = (\va \star \vb) \scalarprod X.
\end{equation}
\end{lemma}
\prf
By associativity of convolution, we have
\begin{equation}
\va \scalarprod (\vb \scalarprod X) 
	= \va \star (\vb \star X) = (\va \star \vb) \star X = (\va \star \vb) \scalarprod X,
\end{equation}
which proves the claim.
\qed

\begin{lemma} \label{lemma:triple-associative}
For any tensor $T \in \Tensor$, vector $\va \in \Vector$ and matrix $X \in \Matrix$, we have
\begin{equation}
T \star (\va \scalarprod X) = \va \scalarprod (T \star X).
\end{equation}
\end{lemma}
\prf
We use associativity and commutativity properties of convolution:
\begin{eqnarray}
T \star (\va \scalarprod X)
	& = & T \star (\va \star X) \\
	& = & (T \star \va) \star X, \ \mbox{ by associativity } \\
	& = & (\va \star T) \star X, \ \mbox{ by {\em commutativity} } \\
	& = & \va \star (T \star X), \ \mbox{ by associativity again } \\
	& = & \va \scalarprod (T \star X),
\end{eqnarray}
which proves the claim.
\qed\\

A natural basis for the group ring $\Matrix = \Vector\GG$ is 
$\widetilde{\basis} = \{\widetilde{B}_{g} = 1_{\Vector}\ket{g} : g \in \GG\}$,
which in Dirac's notation is $\widetilde{B}_{g} = \ket{1_{\GG}}\bra{g}$. But to show that $\Matrix$ is 
a free module under the scalar product $\va \scalarprod X$, Braman \cite{b08} showed that a {\em transposed} 
basis for $\Matrix$ is more suitable.

\begin{figure}[t]
\begin{center}
\epsfig{file=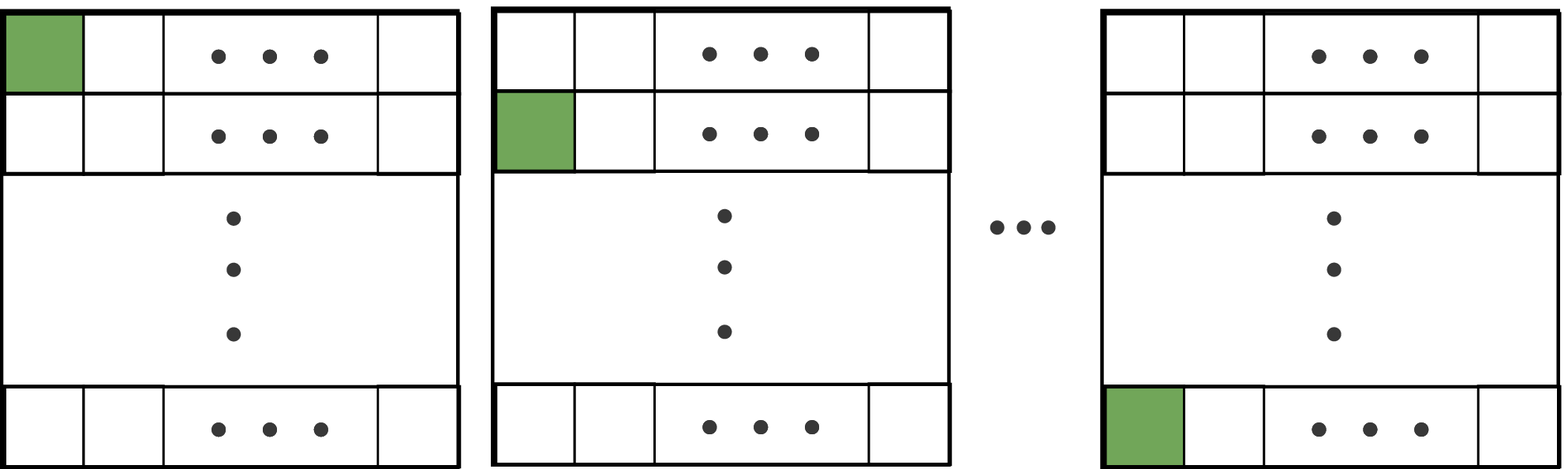, scale=0.5}
\caption{The {\em transposed} basis $\basis = \{B_{g} = \ket{g}\bra{0} : g \in \Zn\}$ 
for the module $\Matrix = \Vector\Zn$ with scalar product $\va \scalarprod X$ (used by Braman \cite{b08}). 
A more natural but unsuitable basis is $\{\widetilde{B}_{g} = \ket{0}\bra{g} : g \in \Zn\}$.
}
\label{figure:matrix_basis}
\end{center}
\end{figure}

\begin{theorem} \label{thm:free-module}
The group ring $\Matrix$ together with the scalar product $\va \scalarprod X$ is a free unitary (left) 
$\Vector$-module, where the basis is $\basis = \{B_{g} = \ket{g}\bra{1_{\GG}} : g \in \GG\}$.
\end{theorem}
\prf
We verify the unitary module properties of $\Matrix = \Vector\GG$ with the scalar product $\va \circ X = \va \star X$.
\begin{enumerate}
\item Scalar product distributes over vector addition: 
\begin{equation}
\va \scalarprod (X + Y) = \va \scalarprod X + \va \scalarprod Y.
\end{equation}
This follows by distributive law for convolution.

\item Scalar product distributes over scalar addition: 
\begin{equation}
(\va + \vb) \scalarprod X = \va \scalarprod X + \vb \scalarprod X.
\end{equation}
This again follows by distributive law for convolution.

\item Scalar product satisfies: 
\begin{equation}
\va \scalarprod (\vb \scalarprod X) = (\va \star \vb) \scalarprod X.
\end{equation}
This follows from Lemma \ref{lemma:scalar-assoc}.

\item Scalar product possesses an identity element:
\begin{equation}
1_{\Vector} \scalarprod X = 1_{\Vector} \star X = X.
\end{equation}
This follows since $1_{\Vector}$ is the identity element of the group ring $\Vector = \Ring\GG$.

\end{enumerate}
To show that $\Matrix$ is a free module, we show that $\basis$ is a basis for $\Matrix$.
We prove that for any matrix $X \in \Matrix$, there are vectors $\va_{\ell} \in \Vector$ so that
$X = \sum_{\ell} \va_{\ell} \scalarprod B_{\ell}$. 
We have
\begin{eqnarray}
\bra{j}X\ket{k} 
	& = & \bra{j} \sum_{\ell} (\va_{\ell} \star B_{\ell})\ket{k},
		\ \ \mbox{ by definition of $\va_{\ell} \scalarprod B_{\ell}$ } \\
	& = & \bra{j} \sum_{\ell} \sum_{i} (\va_{\ell})_{ki^{-1}} B_{\ell}\ket{i},
		\ \ \mbox{ by definition of convolution } \\
	\label{eqn:free-basis}
	& = & \bra{j} \sum_{\ell} \sum_{i} (\va_{\ell})_{ki^{-1}} \ket{\ell}\braket{1_{\GG}}{i},
		\ \ \mbox{ by definition of $B_{\ell}$ } \\
	& = & (\va_{j})_{k}.
\end{eqnarray}
Thus\footnote{Using the natural basis $\widetilde{B}_{\ell} = \ket{1_{\GG}}\bra{\ell}$ 
in Equation (\ref{eqn:free-basis}), yields the degeneracy 
$\bra{j}X\ket{k} = (\sum_{\ell} (\va_{\ell})_{k\ell^{-1}}) \braket{j}{1_{\GG}}$. 
This allows only the first row of $X$ to be defined.},
define $\va_{\ell}$ so that $(\va_{\ell})_{k} = \bra{\ell}X\ket{k}$.
\qed

\section{Isomorphism} \label{section:isomorphism}

In this section, we show that the set of tensors $\Tensor$ is isormophic to the set of module homomorphisms 
over the module of matrices $\Matrix$. As proved by Braman \cite{b08}, this generalizes the fundamental
connection between the set of $n \times n$ matrices and the set of linear transformations over an $n$-dimensional
vector space.

\begin{theorem}
Let~ $\GG$ be an abelian group and~ $\Ring$ be a commutative ring with identity.
Also, let~ $\Vector = \Ring\GG$ be a group ring and~ $\Matrix = \Vector\GG$ be a free module with the 
scalar product $\va \scalarprod X = \va \star X$. 
If~ $\Tensor = M_{n}(\Ring)\GG$ is the group matrix ring, where $M_{n}(\Ring) \cong \hom_{\Ring}(\Vector)$, then
\begin{equation}
\Tensor \ \cong \ \hom_{\Vector}(\Matrix).
\end{equation}
\end{theorem}
\prf
First, we show any tensor $T \in \hom_{\Ring}(\Vector)\GG$ is a module homomorphism over $\Matrix = \Vector\GG$.
Let $A,B \in \Vector\GG$ and $v \in \Vector$. 
Then, $T \star (A+B) = T \star A + T \star B$ holds since convolution distributes over addition.
Also, $T \star (v \scalarprod A) = v \scalarprod (T \star A)$ holds by Lemma \ref{lemma:triple-associative}.
This shows $T \in \hom_{\Vector}(\Vector\GG)$.

Next, we show any module homomorphism $L \in \hom_{\Vector}(\Matrix)$ has a representation as a
tensor $T \in \Tensor$. 
By Theorem \ref{thm:free-module}, $\{B_{g} = \ket{g}\bra{1_{\GG}} : g \in \GG\}$ is a basis for the module $\Matrix$.
Now, let $L \in \hom_{\Vector}(\Vector\GG)$ be so that its action on the basis element $B_{h}$ is defined by
\begin{equation} \label{eqn:lbasis}
L(B_{h}) = \sum_{g \in \GG} \alpha_{g,h} \scalarprod B_{g}, 
\ \ \ \mbox{ where $\alpha_{g,h} \in \Vector$. }
\end{equation}
Then, for any $A \in \Matrix$ with $A = \sum_{h \in \GG} a_{h} \scalarprod B_{h}$, we have
\begin{equation}
L(A) 
= L(\sum_{h \in \GG} a_{h} \scalarprod B_{h})
= \sum_{h \in \GG} L(a_{h} \scalarprod B_{h})
= \sum_{h \in \GG} a_{h} \scalarprod L(B_{h}),
\end{equation}
since $L$ is a module homomorphism over $\Matrix = \Vector\GG$.
Therefore,
\begin{eqnarray}
L(B_{h})\ket{k}
	 & = & \sum_{j \in \GG} \alpha_{j,h} \star B_{j}\ket{k},
	 	\ \ \mbox{ by equation (\ref{eqn:lbasis}) } \\
	 \label{eqn:basis-iso}
	 & = & \sum_{j \in \GG} \sum_{i \in \GG} (\alpha_{j,h})_{ki^{-1}} B_{j}\ket{i},
	 	\ \ \mbox{ by definition of convolution } \\
	 & = & \sum_{j \in \GG} (\alpha_{j,h})_{k} \ket{j},
	 	\ \ \mbox{ since $B_{j} = \ket{j}\bra{1_{\GG}}$. }
\end{eqnarray}
So\footnote{Using the natural basis $\widetilde{B}_{j} = \ket{1_{\GG}}\bra{j}$ in Equation (\ref{eqn:basis-iso}), 
leads to the degeneracy $L(\widetilde{B}_{h})\ket{k} = (\sum_{j} (\alpha_{j,h})_{kj^{-1}})\ket{1_{\GG}}$.},
define a tensor $T \in \hom_{\Ring}(\Vector)\GG$, where $T = \sum_{k \in \GG} T_{k} \ket{k}$ so that
the action of the homomorphism $T_{k}$ on the basis element $\ket{h}$ of $\Vector$ is given by
\begin{equation}
T_{k}(\ket{h}) = \sum_{j} (\alpha_{j,h})_{k} \ket{j}.
\end{equation}
So, $L(B_{h})$ and $T \star B_{h}$ are the same homomorphism (or matrix), since for each $k$:
\begin{equation}
(T \star B_{h})\ket{k} 
	= \sum_{\ell \in \GG} T_{k\ell^{-1}}(B_{h}\ket{\ell}) 
	= \sum_{\ell \in \GG} T_{k\ell^{-1}}(\ket{h}\braket{1_{\GG}}{\ell})
	= T_{k}(\ket{h})
\end{equation}
Therefore, using Lemma \ref{lemma:triple-associative}, we have
\begin{equation}
T \star A 
	= \sum_{h \in \GG} T \star (a_{h} \scalarprod B_{h}) 
	= \sum_{h \in \GG} a_{h} \scalarprod (T \star B_{h}), 
	= \sum_{h \in \GG} a_{h} \scalarprod L(B_{h})
	= L(A).
\end{equation}
This completes the isomorphism between $\hom_{\Ring}(\Vector)\GG$ and $\hom_{\Vector}(\Vector\GG)$.
\qed

\section{Diagonalization} \label{section:diagonalization}

This section shows that the {\em eigenvalues} and {\em eigenvectors} of a tensor can be recovered
provided we are given a diagonalizing equation for the tensor. Although the theorem does not provide a method 
to diagonalize a tensor, it shows that there is a familiar underlying eigenvalue-eigenvector equations that can 
be extracted from the diagonalizing equation.

A module-homomorphism $L \in \hom_{\Ring}(\Vector)$ is called {\em diagonal} if 
for all basis elements $\ket{k}$ of the group ring $\Vector$, 
we have $L(\ket{k}) = \ell_{k} \star \ket{k}$, for some $\ell_{k} \in \Ring$.
A tensor $\td \in \hom_{\Ring}(\Vector)\GG$ is called {\em diagonal} if 
$\td = \sum_{g \in \GG} D_{g} \ket{g}$, 
with each $D_{g}$ being diagonal. 
We use juxtaposition to denote composition of maps.

\begin{figure}[t]
\begin{center}
\epsfig{file=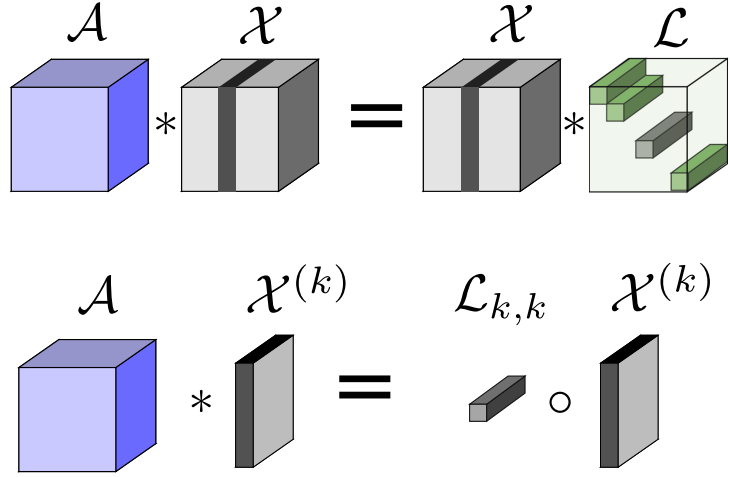, scale=0.75}
\caption{Tensor diagonalization and the eigenvalue-eigenvector equations.}
\label{figure:diagonalization}
\end{center}
\end{figure}

\begin{theorem}
Let~ $\ttt = \sum_{g} T_{g}\ket{g}$, 
$\tx = \sum_{g} X_{g}\ket{g}$, and 
$\tl = \sum_{g} L_{g}\ket{g}$ be tensors in $\Tensor$.
Also, assume that $\tl$ is diagonal with $L_{g}(\ket{k}) = r_{g} \star \ket{k}$.
For~ $k \in \GG$, let $X^{(k)} \in \Matrix$ and $\tl_{k,k} \in \Vector$ be the following formal sums:
\begin{equation}
X^{(k)} = \sum_{g \in \GG} X_{g}(\ket{k}) \ket{g},
\hspace{.5in}
\tl_{k,k} = \sum_{g \in \GG} r_{g} \ket{g}, 
\ \ 
\mbox{ $r_{g} \in \Ring$}.
\end{equation}
If~ $\ttt \star \tx = \tx \star \tl$, 
~then 
\begin{equation}
\ttt \star \tx^{(k)} 
= \tl_{k,k} \scalarprod \tx^{(k)},
\end{equation}
for each $k \in \GG$.
\end{theorem}
\prf
Note for any $h,k \in \GG$, we have 
\begin{equation}
(\ttt \star \tx)_{h}(\ket{k}) 
	= \sum_{g \in \GG} T_{hg^{-1}} \tx_{g}(\ket{k}) \\
	= \sum_{g \in \GG} T_{hg^{-1}} \tx^{(k)}_{g} 
	= (\ttt \star \tx^{(k)})_{h}.
\end{equation}
On the other hand,
\begin{eqnarray}
(\tx \star \tl)_{h}(\ket{k}) 
	& = & \sum_{g \in \GG} X_{hg^{-1}} L_{g}(\ket{k}), \ \mbox{ by convolution in $\Tensor$ } \\
	& = & \sum_{g \in \GG} X_{hg^{-1}}(r_{g} \star \ket{k}), \ \mbox{ by definition of $L_{g}$ } \\ 
	& = & \sum_{g \in \GG} r_{g} \star X_{hg^{-1}}(\ket{k}), \ \mbox{ since $X_{hg^{-1}}$ is a homomorphism } \\
	& = & \sum_{g \in \GG} r_{g} \star \tx^{(k)}_{hg^{-1}}, \ \mbox{ by definition of $\tx^{(k)}$ } \\
	& = & (\tl_{k,k} \star \tx^{(k)})_{h}, \ \mbox{ by definition of $\tl_{k,k}$. }
\end{eqnarray}
This shows that $\ttt \star \tx^{(k)} = \tl_{k,k} \scalarprod \tx^{(k)}$.
\qed

\section{Conclusions} \label{section:conclusions}

In this work, we generalized Braman's work \cite{b08} where it was proved that third-order tensors are exactly
the set of linear transformations over a module formed by second-order tensors (matrices) with scalars that are 
first-order tensors (vectors). Our generalization is based on viewing Braman's construction in an
algebraic framework involving group rings. An advantage of this viewpoint is that the different circulant
products used in \cite{b08} arise naturally as multiplications in various group rings. This is because
convolution is the natural multiplication in these group rings and a circulant-based tensor product is
a convolution based on the cyclic group $\Zn$.

In our generalization, we are able to allow any abelian group $\GG$ to replace the cyclic group $\Zn$
(which defined the convolution), any commutative ring with identity to replace the field $\Real$ of real numbers, 
and arbitrarily high-order tensors to replace third-order tensors provided they form a {\em commutative} ring.
Within this algebraic setting, we also identified Braman's clever choice of the basis for the module of matrices,
which is not the natural basis for the group ring.

\section*{Acknowledgments}

Carmeliza Navasca is supported by National Science Foundation grant DMS-0915100.

\end{document}